# On the ergodicity assumption in Performance-Based Engineering

André T. Beck, Rubia M. Bosse, Isabela D. Rodrigues

# On the Ergodicity Assumption in Performance-Based Engineering

**André T. Beck[1]**
Department of Structural Engineering, University of São Paulo
Av. Trabalhador São-carlense, 400, 13566-590 São Carlos, SP, Brazil. atbeck@sc.usp.br

**Rubia Mara Bosse**
Department of Structural Engineering, University of São Paulo
Av. Trabalhador São-carlense, 400, 13566-590 São Carlos, SP, Brazil. rubiabosse@usp.br

**Isabela Durci Rodrigues**
Department of Structural Engineering, University of São Paulo
Av. Trabalhador São-carlense, 400, 13566-590 São Carlos, SP, Brazil. idrodrigues@usp.br

**Abstract**

In the Performance-Based Engineering (PBE) framework, uncertainties in system parameters, or modelling uncertainties, have been shown to have significant effects on capacity fragilities and annual collapse rates of buildings. Yet, since modelling uncertainties are non-ergodic variables, their consideration in failure rate calculations offends the Poisson assumption of independent crossings. This problem has been addressed in the literature, and errors found negligible for small annual collapse failure rates. However, the errors could be significant for serviceability limit states, and when failure rates are integrated in time, to provide lifetime failure probabilities. Herein, we present a novel formulation to fully avoid the error in integration of non-ergodic variables. The proposed product-of-lognormals formulation is fully compatible with popular fragility modelling approaches in PBE context. Moreover, we address collapse limit states of realistic reinforced concrete buildings, and find errors of the order of 5 to 8% for 50-year lifetimes, up to 14% for 100 years. Computation of accurate lifetime failure probabilities in a PBE context is clearly important, as it allows comparison with lifetime target reliability values for other structural analysis formulations.

---

[1] Corresponding author

1. **Introduction: time-variant reliability analysis**

The theoretical background of Performance-Based Engineering (PBE) is given by time-variant reliability theory. As PBE becomes mainstream and is integrated into design practice, it is important to recall some of its theoretical limitations. One such limitation relates to the integration of failure rates in time, to obtain lifetime probabilities of failure.

Structural systems like bridges and buildings are subject to environmental loads of large magnitude but low probability of occurrence, such as earthquakes and windstorms. The probability that these structures survive application of such loads can be cast as a time-variant reliability problem [1]. This is the probability that a stochastic load process (or load effect) $S(t)$ exceeds the structural capacity within the lifetime $t_D$ of the structure. Loss of equilibrium failure is characterized at the first time the load process crosses the capacity level from below, or first overload. If $\eta(r,t)$ is the arrival rate of the first up-crossing, and if such crossings are independent, the lifetime failure probability $p_f(r,t_D)$ is:

$$p_f(r,t_D) = 1 - \mathcal{R}(r,t_D) = 1 - \mathcal{R}_0 P[N^+(r,t_D) = 0]$$
$$= 1 - (1 - p_{f0})\left(1 - exp\left[-\int_0^{t_D} \eta(r,t)dt\right]\right)$$
$$= p_{f0} + (1 - p_{f0})\left(1 - exp\left[-\int_0^{t_D} \eta(r,t)dt\right]\right) \quad (1)$$

where $\mathcal{R}$ is the reliability, $\mathcal{R}_0$ is the reliability at time $t = 0$, $N^+(r,t_D)$ is the number of capacity crossings within life $t_D$, $p_{f0}$ is the probability that the load process starts in the failure domain, and $r$ is the vector of structural parameters. Eq. (1) describes a Poisson counting process and is valid if the barrier crossings are independent events. Rate $\eta(r,t)$ is time-variant to account for non-stationary load processes and capacity degradation in time (non-homogeneous Poisson process). In most practical applications, the initial failure condition can be disregarded, as $p_{f0} \ll 1$. In this case, Eq. (1) becomes:

$$p_f(r,t_D) = 1 - exp\left[-\int_0^{t_D} \eta(r,t)dt\right] \quad (2)$$

Eqs. (1) and (2) are conditional on a given capacity, or realization $R = r$ of the vector of random system parameters $R$. System parameters $R$ are non-ergodic variables, because they "do not renew at each earthquake event", citing [2]. In a PBE context, uncertain system parameters are mass, damping; the parameters which characterize the load-displacement curve: stiffness, deformation capacity, ultimate strength, and post-peak softening response; as well as model uncertainties of structural response, including parameters of capacity fragilities. In contrast, all variables which are characterized by record-to-record variability, such as load intensity, arrival rates, ground motion time-histories and frequency content, source, path, and site effects, are ergodic in nature, and supposed to be included in crossing rate $\eta(\mathbf{r},t)$.

The unconditional failure probability for lifetime $t_D$, $p_f(t_D)$, is given by the total probability theorem [1]:

$$p_f(t_D) = E_{\mathbf{R}}[p_f(\mathbf{r},t_D)] = \int_{\mathbf{R}} \left(1 - exp\left[-\int_0^{t_D} \eta(\mathbf{r},t)\mathrm{d}t\right]\right) f_{\mathbf{R}}(\mathbf{r})\mathrm{d}\mathbf{r} \qquad (3)$$

where $E_{\mathbf{R}}[\,\cdot\,]$ is the expected value operator, and $f_{\mathbf{R}}(\mathbf{r})$ is the joint probability density function of system variables.

If one manipulates Eq. (1) to isolate $\eta(\mathbf{r},t)$, it becomes clear [3], [4] that $\eta(\mathbf{r},t)$ is a conditional failure rate: the process should start below barrier level $y = f(\mathbf{r})$ at time $t = 0$ and there should be no up-crossings before time $t$. These conditions are important to understand the limitations of the PBE framework, as follows. When we replace arrival rate $\eta(\mathbf{r},t)$ by the rate $\nu^+(y,t)$ at which a stochastic process $S(t)$ crosses barrier level $y$, or the annual rate of failure $\lambda$ due to arrival of earthquakes or windstorms of random magnitude, we may offend the Poisson assumption of independent crossings. There are typically three situations in which this may occur: a) low barrier levels, for which the initial condition is not respected with non-negligible probability; b) clumps of crossings by narrowband load processes; c) random barrier levels. Cases a) and b) above have been addressed in a classic Vanmarcke paper [5], which provides corrections to make $\nu^+(y,t)$ closer to $\eta(y,t)$. Following Vanmarcke [5], Eqs. (1) and (2) are asymptotically exact as barrier level $y \to \infty$. Case c) above has been addressed in [6], [7], as detailed later on in this paper. Case c) addresses the effects of non-ergodic system parameters, which have significant impact on annual collapse rates for buildings subject to earthquakes [8], [9], [10], [11], [12], [13], [14].

2. The ensemble failure rate error in Performance-Based Engineering

The PBE framework is a practical way of implementing time-variant reliability analysis in codified design of structures subject to load processes of extreme magnitude but low probabilities of occurrence. One of the practical features of PBE is the separation of hazard, structural, damage and loss analyses, which can be performed independently by specialized groups. Another relevant feature is consideration of all levels of hazard intensities, by means of the intensity measure (IM) variable (first mode accelerations for earthquakes, wind velocities, etc.). The formulation in Section 1 becomes closer to the PBE formulation by making Eqs. (1) to (3) also dependent on the intensity measure. Eq. (3) in particular becomes:

$$p_f(t_D) = E_R\left[1 - exp\left(-\int_0^{t_D} E_{RTR}[\eta(r, im, t)]dt\right)\right] \quad (4)$$

where $E_{RTR}[\cdot]$ is the expected value over all parameters with record-to-record (RTR) variability (intensity, time-history, spectral shape). In particular, by replacing $\eta(r, im, t)$ by an annual failure rate $\lambda(r, im, t)$, conditional on the intensity measure realization $IM = im$, one obtains:

$$p_f(t_D) = \int_R \left[1 - exp\left(-\int_0^{t_D}\int_{IM} \lambda(r, im, t)h(im)dim\, dt\right)\right] \quad (5)$$

where $h(im) = |dH/dim|$ is the hazard "density" function and $H(im) = P[IM \geq im]$ is the Hazard function. Note that integration over $im$ leads to the annual failure rate:

$$\lambda(r, t) = \int_{IM} \lambda(r, im, t)h(im)dim \quad (6)$$

At this point, one fundamental difference in handling ergodic and non-ergodic variables becomes evident [2], [15], [14]. The expectation over hazard intensity measure (and other variables with RTR variability, for that matter) can be taken inside the time integral in Eq. (5) because these are ergodic variables, whose value changes from one load application to the next. The expectation over non-ergodic variables needs to be taken "outside" the time integral (see Eq. (5)), following the total probability theorem [1].

The expected values in Eqs. (3) to (5) are typically evaluated by Monte Carlo simulation (MCS). When MCS is employed, Eqs. (3) to (5) imply that each independent sample of random system parameters should be exposed to the whole suit of recorded or artificial ground motions. This so-called "robust reliability assessment" [14] involves $N \times N_m$ non-linear structural dynamics time-history computations, where $N$ is the number of ground motion records and $N_m$ is the number of system parameter samples. The so-called "average fragility assessment" [14], versions of which were employed in the literature [8], [9], [10], [12], [16], takes the expectation over $R$ in capacity fragilities, or in annual failure rate calculations:

$$\lambda_E(t) = \int_{IM}\int_R \lambda(r, im, t)dr\, h(im)dim \quad (7)$$

Subscript $(\cdot)_E$ is used in Eq. (7) to stress that this is an *ensemble* failure rate, which cannot be used in Eqs. (1) or (2) without significantly offending the Poisson assumption of independent crossings [2], [14]. Ensemble annual failure rates are also obtained when one-to-one assignment of structural modelling parameters realizations and ground motion records is considered [14]. In this case, one significant advantage is a reduction to $N$ in the total number of dynamic structural analysis.

The non-ergodicity issue has been addressed by Der Kiureghian [2], who stated the error would be negligible for failure rates of the order $\lambda_E \approx 0.01$, up to 20% for $\lambda_E \approx 0.05$, and up to about 30% for $\lambda_E \approx 0.1$. The issue was also addressed by Jalayer and Ebrahimian [14], who found the error would be small for small values of the product of time interval ($t$) by the variance of time-invariant LS probability. The authors also pointed out that the error could be avoided by explicilty considering non-ergodic variables in the format of Eq. (5). The contribution of this manuscript can be seen as a practical way of accomplishing the separation of effects of ergodic and non-ergodic variables, with the purpose of solving Eq. (5), for any lifetime $t_D$.

The non-ergodicity issue was also studied in refs. [3], [6], [7], [17], but in a context of stationary Gaussian processes crossing random Gaussian barriers. In [6], it is shown that the so-called "ensemble crossing rate" error depends not only on the magnitude of the crossing rate, but on the distance between the means of load process and barrier ($\mu_R - \mu_S$), and on the relative magnitude of the standard deviation of the load process ($\sigma_S$), with respect to the barrier ($\sigma_R$). The error is expected to be small when $\mu_R - \mu_S$ is large, and when $\sigma_S \gg \sigma_R$. In a PBE framework, this means that the distance between the mean of intensity measure capacity and the mean of the hazard "density" needs to be large. In earthquake engineering, the impact of random system parameters is significant [8], [9], but we typically have $\sigma_S \gg \sigma_R$. Hence, for problems involving loads of large intensity but low probability of occurrence, the ensemble failure rate error will not be as large as the "order of magnitude" errors found in [6], [7]. Moreover, Beck [7] found that the ensemble crossing rate error was proportional to the number of load cycles. In earthquake engineering, the expected number of load cycles is small, since the actual seismic action contributes a few cycles, which are reduced by a low lifetime probability of occurrence.

The above attenuants by no means imply that the ensemble crossing rate error can be dismissed in PBE. In this manuscript, we address realistic fragility curves for RC frame structures, and show that errors of up to about 10% can still be observed. However, we also show a practical way of avoiding the ensemble failure rate error in PBE. In spite of the above, it is not uncommon to find published papers suggesting that ensemble failure rates can be integrated in time. See for instance [18]–[20], and references cited in [2], paragraph before Eq. (5).

## 3. Consideration of system parameter uncertainty in PBE

Before proceeding, let us recall the SAC-FEMA formulation of PBEE [21]–[23]. Herein, we do not address decision nor damage variables. We evaluate annual limit state (LS) failure rates and lifetime failure probabilities. For a given LS (serviceability, damage, collapse), the annual failure rate is:

$$\lambda_{LS} = \int_0^\infty \left[ \int_0^\infty F(EDP_C | edp, t) \cdot f(edp|im) \mathrm{d}edp \right] h(im) \mathrm{d}im \qquad (8)$$

In Eq. (8), $F(\cdot)$ represents the cumulative distribution function (CDF), $f(\cdot)$ represents the probability density function (PDF), $EDP_C$ is the LS capacity expressed in terms of engineering demand parameters (EDP, interstory drifts, displacements, accelerations, …). The LS capacity is the $edp$ value at which a given limit state (serviceability, damage, collapse) is reached. The intensity measure $im$ may be spectral acceleration or peak ground acceleration. The term inside brackets is the CDF $F_Z(\cdot)$ of the $im$ capacity for a given LS, or fragility function; which provides the probability that a given LS is reached, conditional on $im$. Hence, Eq. (8) can be simplified to:

$$\lambda_{LS} = \int_0^\infty F_Z(IM_C|im) h(im) \mathrm{d}im \qquad (9)$$

The intensity measure capacity $IM_C = Z$ is a function of both ergodic and non-ergodic variables. If non-ergodic variables are included in the rate calculation, then Eq. (9) is equivalent to Eq. (7) [14].

It is convenient to write $Z$ as the product of ergodic ($X$) and non-ergodic ($Y$) variables, as done in [8], [9], [24]:

$$Z = g(X, Y) = X \cdot Y \qquad (10)$$

Here, $X$ represents the effect of all ergodic variables, but the hazard intensity. These variables lead to record-to-record variabilities. Variable $Y$ is the effect, in terms of intensity measure capacity, of random system parameters: $Y = w(\mathbf{R})$. Eq. (10) is convenient because fragilities are typically approximated as lognormal CDF´s. Assuming that $X$ and $Y$ can be approximated as lognormal variables, we have a product of two lognormals, which is also lognormal; and since $\ln(Z)$ has normal distribution, we have:

$$E[\ln(Z)] = E[\ln(X)] + E[\ln(Y)] = \widehat{m}_X + \widehat{m}_Y$$

$$Var[\ln(Z)] = Var[\ln(X)] + Var[\ln(Y)] = \sigma_{lnX}^2 + \sigma_{lnY}^2 \qquad (11)$$

Eq. (10) was employed by [8] and [9] to study the impact of uncertainties in system parameters (called "modelling" uncertainties in these references) on annual failure rates (Eq. 7). The authors show that $\sigma_{lnX}^2$ and $\sigma_{lnY}^2$ are of the same magnitude, and that system parameter uncertainties have a large impact on annual collapse failure rates. The authors proposed two practical ways of evaluating the effects of uncertain system parameters in fragility responses and annual failure rates. Herein, we show a simple procedure for avoiding the ensemble failure rate error, which is compatible with their formulations.

### 4. Avoiding the ensemble crossing rate error

From Eq. (10) one can write: $x = g^{-1}(z, y) = z/y$ and $\partial g^{-1}/\partial z = 1/y$. The CDF of $Z$ can be obtained as [25]:

$$F_Z(z) = \int_0^\infty \int_0^{g^{-1}} f_{XY}(x, y) \mathrm{d}x\, \mathrm{d}y = \int_0^\infty \int_0^z \frac{1}{y} f_{XY}\left(\frac{z}{y}, y\right) \mathrm{d}z\, \mathrm{d}y \tag{12}$$

The effect of record-to-record seismic variables $X$ can be assumed independent of the effect of system variables $Y$ (this is implied in Eqs. (3) to (5)). Hence:

$$F_Z(z) = \int_0^\infty \int_0^z \frac{1}{y} f_X\left(\frac{z}{y}\right) f_Y(y) \mathrm{d}z\, \mathrm{d}y \tag{13}$$

Taking the derivative w.r.t. $z$, one obtains:

$$f_Z(z) = \int_0^\infty \frac{1}{y} f_X\left(\frac{z}{y}\right) f_Y(y)\, \mathrm{d}y \tag{14}$$

Now, if we consider the conditional variable $Z(y) = X \cdot y$, from Eq. (10), Eq. (9) becomes:

$$\lambda_{LS}(y) = \int_0^\infty F_Z(IM_C|im, y) h(im) \mathrm{d}im \tag{15}$$

and Eqs. (3) or (4) can be employed without offending the Poisson assumption of independent crossings:

$$p_f(t_D) = \int_0^\infty \left[1 - exp\left(-\int_0^{t_D} \lambda_{LS}(y) \mathrm{d}t\right)\right] f_y(y) \mathrm{d}y \tag{16}$$

Equation (16) combined with Eq. (10) is the main result of this paper, as it allows one to compute accurate lifetime failure probabilities, considering the effect of system uncertainties in a practical way. Equation (16) provides exact results for the lifetime failure probability, if the contribution of ergodic and non-ergodic variables to fragility functions can be "decomposed" following Eq. (10). The above result is by no means limited to the product "decomposition" of ergodic and non-ergodic variables. Any other functional form which allows one to derive $Z(y)$ can be employed. Eq. (10), however, is the most practical way of combining two lognormal distributions. Surely, the above "decomposition" also assumes neglecting any non-ergodicity in seismic hazard assessment [26].

In principle, Eq. (16) can be employed with solutions 1, 2, and 4 (ie, FA, RF, and AFA) described in ref. [14]. Theoretically, the correction proposed by use of Eqs. (10) and (16) is a practical way of accomplishing the correction of AFA results, by post-processing, to achieve agreement with the robust reliability (RR) solution proposed by Jalayer and Ebrahimian [14].

Note that, in contrast to Eqs. (3) and (4), Eq. (16) involves integration of a scalar variable $Y$, obtained by "decomposing" the intensity measure capacity $Z$ in the contribution of ergodic ($X$) and non-ergodic ($Y$) variables. This is exemplified in the sequence using realistic fragility functions obtained by [8] and [9].

Before proceeding, four observations:

- ✓ First, instead of using Eq. (13), it is easier to compute:

$$\lambda_{LS}(y) = \int_0^\infty f_Z(IM_C|im, y)H(im)\mathrm{d}im \quad (17)$$

since Eq. (17) involves a single integral, whereas Eq. (13) involves a double integral. Eq. (17) is a natural extension of a classic reliability result (see Eqs. 1.18 and 1.19 of ref. [1]).

- ✓ Second, Eqs. (16) and (17) can also be employed to calculate exact annual failure rates. To this end, Eq. (17) is computed with $t_D = 1$, resulting in $p_f(t_D = 1)$; and the exact annual failure rate is back-calculated from Eq. (2). The resulting annual failure rate can then be compared with the ensemble rate from Eq. (9).

- ✓ Third, the formulation in Eqs. (1) to (7) accounts for a failure rate which may vary in time; typically, by effects of strength degradation over lifetime. Assume one obtains fragility functions for the intact structure, and fragility functions corresponding to a degraded (corroded or cracked) structure. If one can find a single multiplying variable $Y$ which, combined with a "deterministic" degrading function (as a parameterized stochastic process: $g(t) \cdot Y$), then the formulation presented herein can be employed to evaluate correct lifetime failure probabilities for degrading structures as well. This is left as a suggestion for further research.

## 5. Behavior of ensemble crossing rate error in toy example

In this section, a toy problem example is addressed, with the main objective of showing how the ensemble crossing rate error varies with main problem parameters. This helps to understand the results for realistic RC frames, to be addressed in the sequence.

Consider a hazard rate function with lognormal distribution:

$$H(im) = \eta(1 - F_S(im)) \quad (18)$$

where $\eta$ is the arrival rate of load pulses and $S$ is the random load intensity. The load intensity $S$, the ergodic capacity variables $X$ and the non-ergodic capacity variables $Y$ have lognormal distributions: $S \sim LN(\hat{m}_S, \sigma_{lnS}) = LN(\hat{m}_S, 1)$, $X \sim LN(\hat{m}_X, \sigma_{lnX}) = LN(\ln(4), 0.4)$, $Y \sim LN(\hat{m}_Y, \sigma_{lnY}) = LN(\ln(0.85), 0.4)$. Figure 1 illustrates the cumulative distributions for this problem.

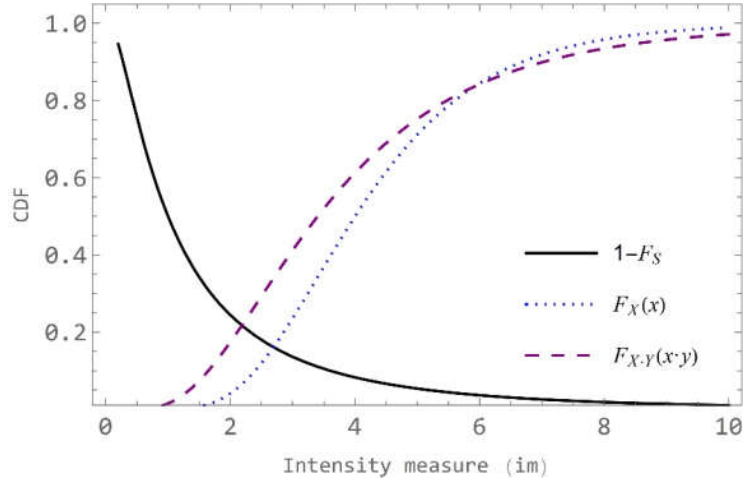

Figure 1: Cumulative distributions of reference toy problem, $\hat{m}_S = \ln(1) = 0$.

The limit state function is conveniently written as:

$$g(s, x, y) = \frac{x\,y}{s} - 1 \tag{19}$$

such that the ensemble failure rate is obtained as:

$$\lambda_E = \eta \Phi(-\beta_{XY})[1 - \Phi(-\beta_{XY})] \quad \text{with} \quad \beta_{XY} = \frac{\hat{m}_X + \hat{m}_Y - \hat{m}_S}{(\sigma_{lnX}^2 + \sigma_{lnY}^2 + \sigma_{lnS}^2)^{1/2}} \tag{20}$$

The conditional failure rate is obtained as:

$$\lambda_{LS}(y) = \eta \Phi(-\beta_X)[1 - \Phi(-\beta_X)] \quad \text{with} \quad \beta_X(y) = \frac{\hat{m}_X + \ln(y) - \hat{m}_S}{(\sigma_{lnX}^2 + \sigma_{lnS}^2)^{1/2}} \tag{21}$$

This last failure rate is integrated over time using Eq. (16).

The percentual error between the ensemble rate and the exact solutions is computed by using the result of Eq. (20) in Eq. (2), with reference to the result of Eq. (16):

$$Err\%p_f(t_D) = 100 \frac{\left((1 - Exp\left[-\int_0^{t_D} \lambda_E dt\right]) - p_f(t_D)\right)}{p_f(t_D)} \tag{22}$$

The percentual error is illustrated in Figure 2, in terms of $\Delta\sigma_{lnY}$ ($\pm 20\%$ variation on standard deviation of non-ergodic variables), for $t_D = 50$. Herein, we are not interested in the absolute value of the error, because problem parameters are arbitrary. Rather, we are interested in the behavior of the error w.r.t. the main problem parameters. As observed in Figure 2, the error grows as $\sigma_{lnY}$ increases, as could be expected. Yet, Figure 2 also shows that the percentual error changes for different hazard rate functions, in terms of median of pulse intensities ($\widehat{m}_S$) and pulse arrival rate ($\eta$). The main conclusion here, with relevance to the following results for realistic structures, is that the trends of ensemble crossing rate errors in terms of capacity parameters should be evaluated for fixed hazard functions.

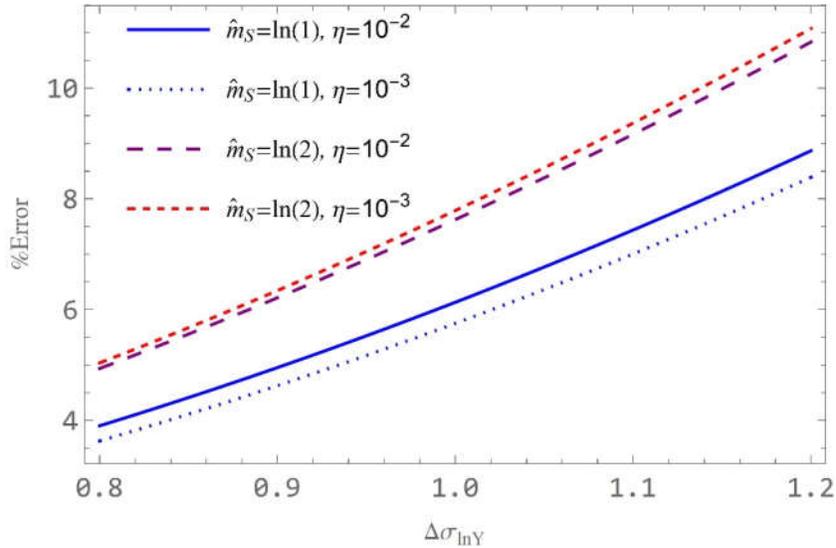

Figure 2: Percentual error in failure probabilities $Err\%p_f$, evaluated with inclusion of non-ergodic variables in failure rate calculation.

The ensemble crossing rate error is believed to be proportional to crossing rates and to failure probabilities since the classical paper of Der Kiureghian [2]. Yet, in ref. [2] only one hazard function was considered. A toy example was addressed considering two values for the variance of non-ergodic capacity variables and several values of the $edp$ parameter. For a fixed hazard function, the ensemble crossing rate error is proportional to $\sigma_{lnY}$ (as shown in Fig. 2) and inversely proportional to the design margin, as shown in [6, 7]. In this situation, the ensemble error is indeed proportional to failure probabilities, as observed in [2]. However, it is shown herein (Figures 3 and 4) that this is not the case when different hazard functions and capacity parameters are considered. Three cases are considered, by changing median pulse intensity ($\widehat{m}_S$) and standard deviation of non-ergodic variables ($\sigma_{lnY}$):

a) $\widehat{m}_S = \ln(1)$ and $\sigma_{lnY} = 0.4$;
b) $\widehat{m}_S = \ln(1)$ and $\sigma_{lnY} = 0.2$;
c) $\widehat{m}_S = \ln(2)$ and $\sigma_{lnY} = 0.2$.

Figure 3 shows lifetime failure probabilities and Figure 4 shows the corresponding percentual ensemble crossing rate errors, for the three cases. Note that failure probabilities increase significantly in time, whereas ensemble crossing errors are nearly constant. In Figure 4 for case c), failure probabilities change by two orders of magnitude (from 0.001 to 0.1), whereas the percentual error is an almost constant 2%.

It is further observed that the relationship between ensemble errors and failure probabilities change with problem parameters (cases a to c above): failure probabilities reduce (marginally) when $\sigma_{lnY}$ is halved (case a to b) and increase when $\hat{m}_S$ is increased (case b to c). Yet, the ensemble errors in Figure 4 show a different trend: they are larger for case a) and significantly smaller for case c). Failure probabilities are larger for case c), but errors are larger for case a). A reduced "design" margin increases failure probabilities (case c), and larger non-ergodic uncertainty increases ensemble error (case a). In time-variant reliability terms, reducing the barrier increases crossing rates and failure probabilities, and increasing barrier uncertainty increases ensemble crossing rate error [6, 7].

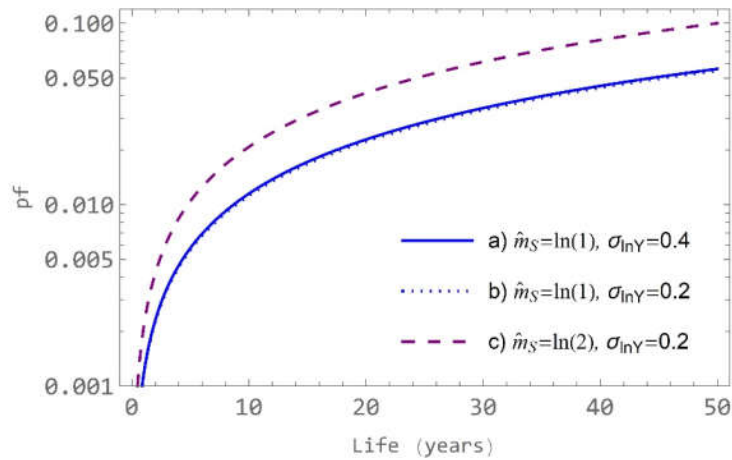

Figure 3: Failure probabilities in time for different parameters of toy problem, $\eta = 10^{-2}$.

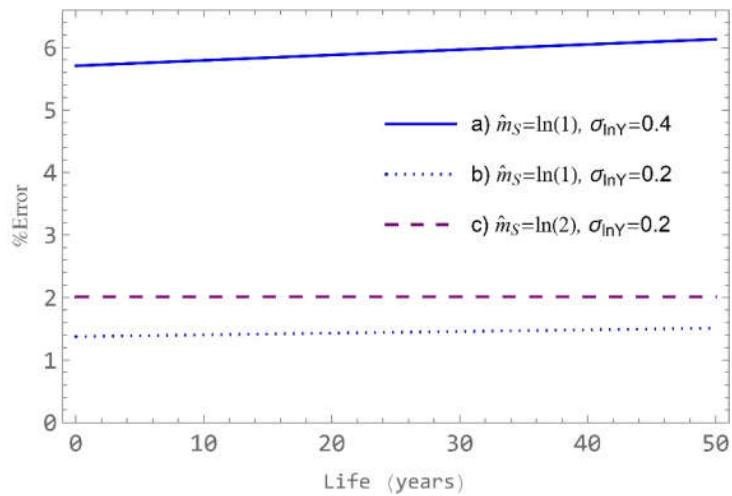

Figure 4: Percentual error in time for different parameters of toy problem, $\eta = 10^{-2}$.

## 6. Realistic fragility curves including system parameter uncertainties

Early studies like [30] addressing effects of uncertainties in damping, mass and material strength on the pre-collapse capacity of reinforced concrete (RC) moment frames found smaller contribution of these variables. Yet, when the collapse capacity is addressed, it is found that the uncertainty in deformation capacity and post-peak softening responses have large impact in conditional collapse fragilities [27].

Haselton & Deierlein [28] evaluated the collapse capacity of 30 ductile RC moment frames between one and twenty stories, located at a near-fault site in Los Angeles. Herein, we address three cases of the 4-storey office building reported in [8]. The relevant data to reconstruct fragility curves is reported in Table 1. This includes the median ($\hat{m}$) and $\sigma_{ln}$ for record-to-record variabilities only ($X$), and including system parameter uncertainties (variable $Z$ in Eq. 10). The median and $\sigma_{ln}$ of non-ergodic ($Y$) variables only is also reported in Table 1, as evaluated from Eq. (11).

Table 1: Data on realistic fragility functions for RC moment frames

| Ref. | Case | Stories | $T_1$ (s) | Margin w.r.t. $S_a(T_1)_{2/50}$ | RTR ($X$) | | Total ($Z$) | | System only ($Y$) | |
|---|---|---|---|---|---|---|---|---|---|---|
| | | | | | $\hat{m}$ | $\sigma_{ln}$ | $\hat{m}$ | $\sigma_{ln}$ | $\hat{m}$ | $\sigma_{ln}$ |
| [8] | 1 | 4 | 1.00 | 2.11 | 1.7 | 0.30 | 1.7 | 0.58 | 1 | 0.5 |
| | 2 | | 1.00 | 2.61 | 2.1 | 0.29 | 2.1 | 0.578 | 1 | 0.5 |
| | 3 | | 1.00 | 3.48 | 2.8 | 0.34 | 2.8 | 0.605 | 1 | 0.5 |
| [9] | 4 | 4 | 1.12 | 1.52 | 1.3 | 0.40 | 1.10 | 0.48 | 0.85 | 0.26533 |
| | 5 | 12 | 2.01 | 1.23 | 0.61 | 0.473 | 0.56 | 0.52 | 0.918 | 0.21603 |
| | 6 | 12 | 1.98 | 0.63 | 0.3 | 0.45 | 0.28 | 0.50 | 0.933 | 0.21795 |
| | 7 | 12 | 2.26 | 0.73 | 0.35 | 0.415 | 0.38 | 0.49 | 1.086 | 0.26053 |

In order to obtain the effects of system parameter uncertainties, Goulet et al. [8] employed a second moment FOSM-based scheme, for which the median of fragility functions do not change, as observed in Table 1. Hence, the only effect of system uncertainties is to increase $\sigma_{lnZ}$, in comparison to $\sigma_{lnX}$. For this case, $\hat{m}_Y = 1$. The three cases reported in Table 1 correspond to: case 1) $\epsilon$-based ground motion selection not considered (see [8] for details); case 2) base case, considering only one horizontal component of ground motions and case 3) using both horizontal components.

Liel [29] studied the collapse capacity of 3 ductile and 3 non-ductile RC frames of one, two, four and twelve stories. Herein, we focus on four of these frames, as detailed in Table 1. Cases 4 and 5 are ductile frames, and cases 6 and 7 are non-ductile, of four and twelve stories, respectively. The site is not-near fault in Los Angeles. Natural/fundamental periods are around 1 and 2 seconds and $Sa(T_1)_{2/50} = 0.46\ g$. Liel et al. [9] presented a response-surface methodology to estimate the effects of system uncertainties in collapse fragilities. As observed in Table 1, this more accurate method produces a shift in medians, as well as an increase in $\sigma_{ln}$. The change in $\sigma_{ln}$ is not as large as that obtained in [8], but combined with reduction in medians, this can have an impact on ensemble failure rate errors. Importantly, this method also resulted in an increase in median for case 7, which is somewhat unexpected.

Vamvatsikos and Fragiadakis [12] observed that the mean fragility response is not the fragility response obtained for mean values of system parameters, a result which is typical of non-linear responses. Jalayer and Ebrahimian [14] studied the three stories RC frame of a school in Italy, and observed a shift in median, and an increase in standard deviations, similar to results in ref. [9]. Fragility curves obtained in ref. [14] are not considered herein because the data is only available in graphical form.

## 7. Results for realistic fragility curves

The methodology proposed in Section 4 is employed to calculate annual failure rates and lifetime failure probabilities for the RC structures described in Section 5. Figure 5 shows the fragility functions for cases (2) and (4) of Table 1. Included in Figure 5 are the collapse capacities for record-to-record variability only, $F_Z(X, \widehat{m}_Y = 1)$; the total collapse capacities for RTR plus system variabilities, $F_{ZE}(X, Y)$, as reported in the references; and the total collapse capacities as calculated from Eqs. (10) and (13), $F_{ZE}(X \cdot Y)$. The last two fragility curves are indicated with an additional subscript $(\cdot)_E$ to stress that these are ensemble CDF´s, already incorporating the uncertainty in system parameters. As observed in Figure 5, Eqs. (10) and (13) provide an exact representation of these fragility curves, when system parameters produce only variance increase (Fig. 5 A) and when they result in mean shift plus variance increase (Fig. 5 B).

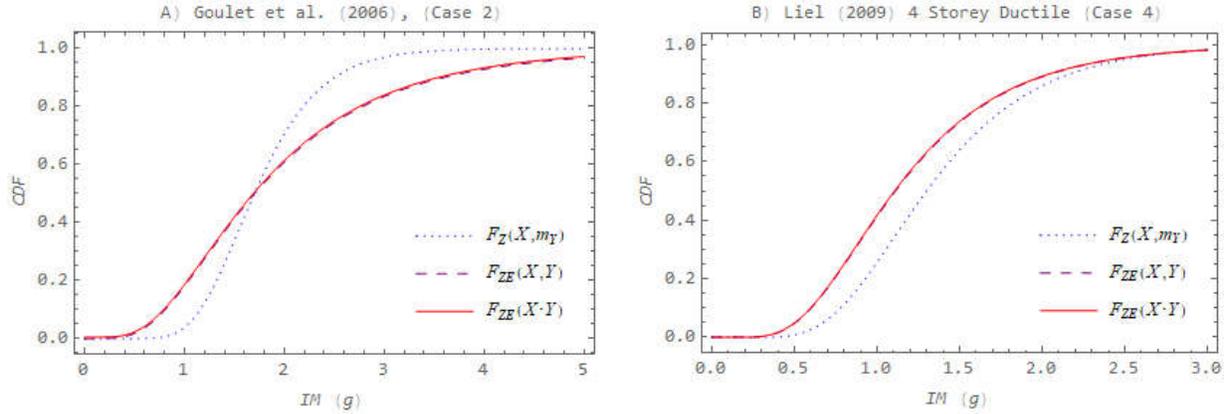

Figure 5: Reproduction of realistic fragility curves for A) Case (2), 4-story frame from [8] and B) Case (4) 4-story frame from [9].

Annual collapse rates for the seven structures listed in Table 1 are shown in Table 2. This includes collapse rates for RTR variability only; exact collapse rates calculated from Eqs. (17) and (2), as explained in Second observation, Section 4; and ensemble collapse rates calculated from Eqs. (7) or (9). Annual collapse rates shown in columns a) and c) are close to the values provided in the references (values are not exactly same due to manual approximation of hazard curves from these references). By dividing the rates in column b) by those in column a), we can appreciate the impact of system uncertainties on annual failure rates. Results obtained by Goulet et al. [8], cases 1 to 3, show the highest impacts, as the logarithmic standard deviations for system variables $Y$ only are larger than for RTR only ($X$) (see Table 1). For the fragility functions obtained by Liel et al. [9], the impact of system uncertainties is smaller. Although Liel et al. [9] observed a reduction in medians (in most cases), the $\sigma_{ln}$ for system variables $Y$ only is smaller than for RTR only ($X$).

The difference between exact and ensemble annual failure rates is minor, as observed in columns b) and c) of Table 2. The percentual error in these results is evaluated from (see Eqs. (7) and (17)):

$$Err\%\lambda = 100 \frac{(\lambda_E - \lambda_{LS}(\hat{m}_y=1))}{\lambda_{LS}(\hat{m}_y=1)} \qquad (23)$$

As observed, the largest error for the structures studied herein is less than 0.2% (Case 1, with small design margin ($2.11g$) and large standard deviation $\sigma_{lnY}$ with respect to $\sigma_{lnX}$). This confirms the observation by Der Kiureghian [2] that the error is "trivial" for small annual failure rates. This also confirms the observation in [14] that the error is small for small time interval ($t_D = 1$ in Table 2).

Table 2: Annual collapse rates $\lambda_{LS} \times 10^4$ (a, b and c), ratio b/a, and % ensemble failure rate error.

| Frame Case | a) RTR only | b) Total, Eq. (17) (exact) | c) Total, Eqs. (7,9) (with ensemble rate error) | Ratio b/a (effect of system uncertainty) | $Err\%\lambda$ 100(c-b)/b |
|---|---|---|---|---|---|
| 1 | 0.215 | 2.564 | 2.569 | 11.92 | 0.17 |
| 2 | 0.049 | 1.087 | 1.088 | 22.20 | 0.11 |
| 3 | 0.013 | 0.422 | 0.422 | 31.55 | 0.07 |
| 4 | 2.413 | 6.743 | 6.746 | 2.80 | 0.04 |
| 5 | 6.790 | 11.086 | 11.089 | 1.63 | 0.03 |
| 6 | 51.168 | 67.903 | 67.969 | 1.33 | 0.10 |
| 7 | 30.937 | 30.565 | 30.592 | 0.99 | 0.09 |

Failure probabilities computed for a lifetime $t_D = 50$ years are shown in Table 3, for the seven RC frames studied herein. This table also shows results for RTR variabilities only, exact results, and results computed from the ensemble crossing rates. The ratios of failure probabilities (Column c) in Table 3) are about the same observed in Table 2, showing that the impact of system uncertainties in annual failure rates or in lifetime failure probabilities is about the same. Yet, the difference between exact and ensemble failure probabilities is larger, and no longer dismissible. The percentual error in failure probabilities is computed from Eq. (22). As observed in Table 3, the ensemble crossing rate error for a lifetime $t_D = 50$ years, and for the realistic RC frames studied herein, can be as large as about 8%. For lifetime $t_D = 100$ years, the errors range from about 3% to 14%, for the seven cases studied herein, as observed in Table 4.

Table 3: Fifty-year failure probabilities, ratio b/a, % ensemble crossing rate error, error parameter of [7], product $(t_D \cdot Var[p_f])$.

| Frame Case | a) RTR only | b) Total, Eq. (16) (exact) | c) Total, Eq. (7) in Eq. (2) (with ensemble rate error) | Ratio b/a (effect of system uncertainty) | $Err\%p_f$ 100(c-b)/b | Error parameter (Eq. 24) | $t_D \cdot Var[p_f]$ (ref. [14]) |
|---|---|---|---|---|---|---|---|
| 1 | 0.0011 | 0.0119 | 0.0128 | 11.03 | 7.62 | 0.54 | 0.07 |
| 2 | 0.0002 | 0.0052 | 0.0054 | 21.05 | 5.27 | 0.44 | 0.02 |
| 3 | 0.0001 | 0.0020 | 0.0021 | 30.60 | 3.05 | 0.31 | 0.01 |
| 4 | 0.0120 | 0.0326 | 0.0332 | 2.72 | 1.79 | 0.38 | 0.05 |
| 5 | 0.0334 | 0.0532 | 0.0539 | 1.59 | 1.40 | 0.46 | 0.07 |
| 6 | 0.2257 | 0.2772 | 0.2881 | 1.23 | 3.94 | 0.88 | 0.70 |
| 7 | 0.1433 | 0.1363 | 0.1418 | 0.95 | 4.07 | 0.79 | 0.42 |

Table 4: Hundred-year failure probabilities, ratio b/a, % ensemble crossing rate error, error parameter of [7], product ($t_D \cdot Var[p_f]$).

| Frame Case | a) RTR only | b) Total, Eq. (16) (exact) | c) Total, Eq. (7) in Eq. (2) (with ensemble rate error) | Ratio b/a (effect of system uncertainty) | $Err\%p_f$ 100(c-b)/b | Error parameter (Eq. 24) | $t_D \cdot Var[p_f]$ (ref. [14]) |
|---|---|---|---|---|---|---|---|
| 1 | 0.0021 | 0.0222 | 0.0254 | 10.34 | 14.16 | 0.54 | 0.40 |
| 2 | 0.0005 | 0.0098 | 0.0108 | 20.11 | 9.91 | 0.44 | 0.14 |
| 3 | 0.0001 | 0.0040 | 0.0042 | 29.77 | 5.83 | 0.31 | 0.04 |
| 4 | 0.0238 | 0.0631 | 0.0652 | 2.64 | 3.46 | 0.38 | 0.36 |
| 5 | 0.0656 | 0.1022 | 0.1050 | 1.56 | 2.69 | 0.46 | 0.45 |
| 6 | 0.4005 | 0.4635 | 0.4932 | 1.16 | 6.40 | 0.88 | 2.67 |
| 7 | 0.2661 | 0.2456 | 0.2636 | 0.92 | 7.33 | 0.79 | 2.19 |

Confirming results of Section 5, it is noted in Tables 2, 3 and 4 that the ensemble crossing rate error is not proportional to annual failure rates, as supposed by Der Kiureghian [2]. Note that the largest % error is obtained for case 1), which has the third smallest failure rate. The smallest % error is obtained for case 5, which has the third largest annual failure rates. These results can be further interpreted based on [6] and [7]. Although these works address stationary Gaussian processes crossing random Gaussian barriers, most relevant tendencies apply to the present PBE problem. The authors identified an ensemble crossing rate error parameter, which in the present problem can grossly be approximated as:

$$Err\_p_f(t_D) \propto \frac{\sigma^2+1}{\mu} \approx \frac{(\sigma_{lnY}/\sigma_{lnX})^2+1}{D_m/\sigma_{lnX}} \qquad (24)$$

where $D_m$ is the design margin w.r.t. $S_a(T_1)_{2/50}$. As shown in Section 5, values of Eq. 24 cannot be compared across different hazard functions, because the hazard curve is not reflected in Eq. (24). Yet, it can be observed that the error parameters, listed in pre-last columns of Tables 3 and 4, are roughly proportional to the ensemble failure rate errors in the same tables. Among the Goulet et al. [8] frames, the largest error parameter value is obtained for Case 1, for which the largest ensemble crossing rate error was obtained. Among the Liel et al. [9] frames with $T_1 \approx 2s$ (Cases 5 to 7), the smallest error parameter value is obtained for Case 5, for which the smallest ensemble crossing rate error was obtained. Lastly, Case 4 of ref. [9] with $T_1 \approx 1s$ is a loner, which cannot be directly compared to the others (single frame with this hazard function).

These results can also interpreted based on the insight from ref. [14]. Following the authors, the ensemble error is proportional to the product $t_D \cdot Var[p_f]$, where $Var[p_f] = E[p_f^2] - E[p_f]^2$ and:

$$E[p_f^2] = \int_0^\infty \left[ 1 - exp\left(-\int_0^{t_D} \lambda_{LS}(y)dt\right)\right]^2 f_y(y)dy \qquad (25)$$

The last columns of Tables 3 and 4 present the calculated values of this product. Again, we observe that the proposed product is proportional to the ensemble error if analyzed within the same hazard function. For cases 1 to 3, largest product is case 1, which yields largest ensemble rate error. Within cases 5 to 7, smallest product is case 5, which yields smallest ensemble rate error. The comparison between cases 6 and 7 breaks the trends, for the error parameters of Eqs. (24) and (25). One possible reason for this break is that for case 7 an "unusual" increase of median was observed as effect of system uncertainties. We note that Eq. (24) is just a hint or first approximation of ensemble error dependency in PBE. More importantly, a procedure has been presented herein to fully avoid the ensemble crossing rate in PBE lifetime failure probability evaluation.

## 8. Concluding remarks

In this manuscript, we addressed time-variant reliability foundations of performance-based engineering (PBE). Specifically, we addressed the integration of system parameter uncertainties in annual failure rates and lifetime probabilities of failure. By addressing realistic fragility functions for RC frame structures, we showed that the error involved in including non-ergodic variables in lifetime failure probability calculations can be of the order of about 10 to 15%, and that it is not proportional to failure rates, as previously thought. We showed that the error increases when the variance of system parameters increases with respect to variance of record-to-record variables; and when the design margin is smaller. More importantly, we showed that the error in failure probabilities and in failure rate calculations can be avoided by decomposing the response fragility in ergodic and non-ergodic terms, and performing the integration over non-ergodic variables "outside" the time-integration of failure rates. Herein, a practical product-of-lognormals approach was adopted in the ergodic × non-ergodic fragility decomposition. Computation of accurate lifetime failure probabilities in a PBE context is clearly important, as it allows comparison with lifetime target reliability values for conventional structural analysis formulations.

**Acknowledgements**

Funding of this research project by Brazilian agencies CAPES (Brazilian Higher Education Council), CNPq (Brazilian National Council for Research, grant n. 309107/2020-2), FAPESP (grant n. 2020/14072-7) and joint FAPESP-ANID (São Paulo State Foundation for Research - Chilean National Agency for Research and Development, grant n. 2019/13080-9) is also acknowledged.